\documentclass[conference]{IEEEtran}

\usepackage{cite} 
\usepackage{amsthm,amsmath, amsfonts, amssymb} %\usepackage{algorithm, algpseudocode}
\newtheorem{theorem}{Theorem}[section]

\usepackage[isolatin]{inputenc}
\usepackage{enumitem}
\usepackage{subfig}
\usepackage{xcolor}
\usepackage{bm,graphicx, enumerate}
\usepackage{verbatim}
\usepackage{algorithm, algorithmic, tabularx}
\usepackage{dsfont}

\newcommand{\Lap}{\ensuremath{\ma{L}}}
\newcommand{\Fou}{\ensuremath{\ma{U}}}
\newcommand{\fou}{\ensuremath{\vec{u}}}
\renewcommand{\leq}{\ensuremath{\leqslant}}

\newcommand{\adjoint}{\ensuremath{{\intercal}}}

\newcommand{\norm}[1]{\ensuremath{\left\| #1\right\|}}

\newcommand{\ma}[1]{\ensuremath{\mathsf{#1}}}
\renewcommand{\vec}[1]{\ensuremath{\bm{#1}}}

\newtheorem{remark}[theorem]{Remark}

\usepackage{bm}

\newsavebox{\ieeealgbox}
\newenvironment{boxedalgorithmic}
  {\begin{lrbox}{\ieeealgbox}
   \begin{minipage}{\dimexpr\columnwidth-2\fboxsep-2\fboxrule}
   \begin{algorithmic}}
  {\end{algorithmic}
   \end{minipage}
   \end{lrbox}\noindent\fbox{\usebox{\ieeealgbox}}}

\begin{document}
\title{Analyzing the Approximation Error of the \\Fast Graph Fourier Transform}

\author{\IEEEauthorblockN{Luc Le Magoarou\IEEEauthorrefmark{1},
Nicolas Tremblay\IEEEauthorrefmark{2} and 
R\'emi Gribonval\IEEEauthorrefmark{3}}
\IEEEauthorblockA{\IEEEauthorrefmark{1}b\raisebox{0.2mm}{\scalebox{0.7}{\textbf{$<>$}}}com, Rennes, France}
\IEEEauthorblockA{\IEEEauthorrefmark{2}Univ. Grenoble Alpes, CNRS, GIPSA-lab, Grenoble, France}
\IEEEauthorblockA{\IEEEauthorrefmark{3}INRIA Rennes Bretagne-Atlantique, Rennes, France}}

\maketitle

\begin{abstract}
The graph Fourier transform (GFT) is in general dense and requires $\mathcal{O}(n^2)$ time to compute and $\mathcal{O}(n^2)$ memory space to store. In this paper, we pursue our previous work on the approximate fast graph Fourier transform (FGFT). The FGFT is computed via a truncated Jacobi algorithm, and is defined as the product of $J$ Givens rotations (very sparse orthogonal matrices). The truncation parameter, $J$, represents a trade-off between precision of the transform and time of computation (and storage space). %for instance, if $J=\mathcal{O}(n^2)$, the FGFT is exactly the GFT, but there is no computational gain compared to the GFT. 
We explore further this trade-off and study, on different types of graphs, how is the approximation error distributed along the spectrum. 
\end{abstract}

\section{Introduction}
Recently, methods have been developed to analyze and process signals defined over the vertices of a graph \cite{Shuman2013,Sandryhaila2013}, instead of over a regular grid, as is classically assumed in discrete signal processing. The starting point of graph signal processing is to define a graph Fourier transform (GFT), via an analogy with classical signal processing. Depending on the preferred analogy, there exists several definitions of GFTs, all of them based on the diagonalisation of a graph operator, may it be the adjacency matrix~\cite{Sandryhaila2013}, the Laplacian matrix~\cite{Shuman2013}, their degree-normalized versions, etc. We refer the reader to~\cite{tremblay_gsp_chapter} for a recent review on different existing definitions of GFTs. 

In this paper, we restrict ourselves to the study of the Laplacian-based Fourier matrix of undirected graphs. In this case, and denoting by $n$ the number of vertices of the graph, the Laplacian $\mathbf{L}\in\mathbb{R}^{n\times n}$ (see the Methods Section for a precise definition) may be diagonalized as:
\begin{equation}
\mathbf{L} = \mathbf{U}\boldsymbol{\Lambda}\mathbf{U}^\top,
\label{eq:defGFT}
\end{equation} 
where $\mathbf{U} \in \mathbb{R}^{n \times n}$ is an orthonormal matrix whose columns are the graph Fourier modes and $\boldsymbol{\Lambda} \in \mathbb{R}^{n \times n}$ is a non-negative diagonal matrix whose diagonal entries correspond to the graph's generalized frequencies. 

The graph Fourier transform of a graph signal $\bm{x}\in\mathbb{R}^n$ is its scalar product with all Fourier modes, i.e., $\tilde{\bm{x}}=\mathbf{U}^\top\bm{x}\in\mathbb{R}^n$. This Fourier matrix $\mathbf{U}$ being in general dense, these $n$ scalar products require $\mathcal{O}(n^2)$ computations.  Nevertheless, in the classical signal processing case, the well-known Fast Fourier Transform (FFT) \cite{CooleyTukey1965} allows to apply the Fourier transform in only $\mathcal{O}(n \log n)$ arithmetic operations. In fact, the FFT is a fast linear algorithm \cite{Morgenstern1975}, which implies that the classical Fourier matrix can be factorized into sparse factors, as discussed in \cite{Lemagoarou2016a}. It is natural to wonder if this kind of factorization can be generalized to the graph Fourier transform, in order to tend towards a ``fast'' graph Fourier transform, scaling in $\mathcal{O}(n\log{n})$, thus mimicking the classical case. 

In~\cite{magoarou_approximate_2017}, we propose such a factorization approach. More precisely, we develop a truncated Jacobi method that amounts to an approximate diagonalization of the Laplacian $\mathbf{L}$, as
\begin{equation}
\mathbf{L} \approx \mathbf{S}_1\dots\mathbf{S}_J\hat{\boldsymbol{\Lambda}}\mathbf{S}_J^T\dots\mathbf{S}_1^T,
\label{eq:defGFFT}
\end{equation} 
where the so-called Givens rotation matrices $\mathbf{S}_1,\dots,\mathbf{S}_J$ are both sparse \emph{and} orthogonal, and $\hat{\boldsymbol{\Lambda}}$ is a diagonal matrix whose diagonal entries are approximations of the graph frequencies. $J$ is typically chosen such as $J=o(n^2)$. The approximate fast graph Fourier transform of $\bm{x}$ reads  $\hat{\mathbf{U}}^\top\bm{x}=\mathbf{S}_J^\top\dots\mathbf{S}_1^\top\bm{x}$. It is approximate because we stop the Jacobi algorithm before its convergence. It is fast because applying each Givens rotation $\mathbf{S}_i$ to a vector requires only 6 elementary operations (4 multiplications and 2 sums): computing $\hat{\mathbf{U}}^\top\bm{x}$ thus requires only $6J$ computations. 

\textbf{Contributions.} In this paper, we build upon our previous work~\cite{magoarou_approximate_2017} and explore more precisely the approximation error due to the truncated Jacobi method: where in the spectrum are the errors localized, and what are the consequences from a graph signal processing point of view.

%In this paper, we decide to fix the number of Givens rotation $J$ to $n\log{n}$ in order to mimick the classical FFT, and study the approximation performances of the fast graph Fourier transform on different types of graphs. 

\section{Methods}
\label{sec:methods}
\subsection{Notations, conventions and definitions}
\label{sec:notations}

{\noindent \bf General notations.} Matrices are denoted by bold upper-case letters: $\mathbf{A}$; vectors  by bold lower-case letters: $\mathbf{a}$; the $i$th column of a matrix $\mathbf{A}$ by $\mathbf{a}_i$; its entry on the $i$th row and $j$th column by $a_{ij}$. 
Sets are denoted by calligraphic symbols: $\mathcal{A}$, and we denote by $\delta_{\mathcal{A}}$ the indicator function of the set $\mathcal{A}$ in the optimization sense ($\delta_{\mathcal{A}}(x) = 0$ if $x \in \mathcal{A}$, $\delta_{\mathcal{A}}(x) = +\infty$ otherwise). The standard vectorization operator is denoted $\text{vec}(\cdot)$. The $\ell_0$-pseudonorm is denoted $\left\Vert\cdot\right\Vert_0$ (it counts the number of non-zero entries), $\left\Vert\cdot\right\Vert_F$ denotes the Frobenius norm, and $\left\Vert\cdot\right\Vert_{2}$ the spectral norm. By abuse of notations, $\|\mathbf{A}\|_{0} = \|\text{vec}(\mathbf{A})\|_{0}$. The identity matrix is denoted $\mathbf{Id}$. % denotes the identity matrix. \\

{\noindent \bf Graph Laplacian.} We consider in this paper undirected weighted graphs, denoted $\mathcal{G} \triangleq \{\mathcal{V},\mathcal{E},\mathbf{W}\}$, where $\mathcal{V}$ represents the set of vertices (otherwise called nodes), $\mathcal{E} \subset \mathcal{V} \times \mathcal{V}$ is the set of edges, and $\mathbf{W}$ is the weighted adjacency matrix of the graph. We denote $n\triangleq |\mathcal{V}|$ the total number of vertices and the adjacency matrix $\mathbf{W} \in \mathbb{R}^{n \times n}$ is symmetric and such that $w_{ij} = w_{ji}$ is non-zero only if $(i,j)\in \mathcal{E}$ and represents the strength of the connection between nodes $i$ and $j$. We define the degree matrix $\mathbf{D} \in \mathbb{R}^{n \times n}$ as a diagonal matrix with $\forall i$, $d_{ii} \triangleq \sum_{j=1}^n w_{ij}$, and the combinatorial Laplacian matrix $\mathbf{L}\triangleq \mathbf{D} - \mathbf{W}$ (we only consider this type of Laplacian matrix in this paper, and hereafter simply call it the Laplacian).
\\

{\noindent \bf Givens rotations. }  An $n$-dimensional Givens rotation~\cite{Givens1958} is a linear transformation that does not act on $n-2$ coordinates and rotates the two remaining by an angle $\theta\in [0;2\pi[$. Noting $p$ and $q$ the indices of the two rotated coordinates, Givens rotations thus correspond to matrices of the following form,

%\begin{equation*}
%\left( 
%\begin{array}{ccccc}
%\ddots &&&&\\
%&c&&-s&\\
%&&\ddots&&\\
%&s&&c&\\
%&&&&\ddots\\
%\end{array}
%\right),
%\end{equation*}
%
%\begin{equation*}
%\left( 
%\begin{array}{ccccccccccc}
%1&&&&&\\
%&\ddots &&&&\\
%& &1&&&\\
%&&&c&&&&-s&\\
%&&&&1&&\\
%&&&&&\ddots&&\\
%&&&&&&1&\\
%&&&s&&&&c&\\
%&&&&&&&&1\\
%&&&&&&&&&\ddots\\
%&&&&&&&&&&1\\
%\end{array}
%\right),
%\end{equation*}

\begin{figure}[!h]
\centering
\includegraphics[width=0.7\columnwidth]{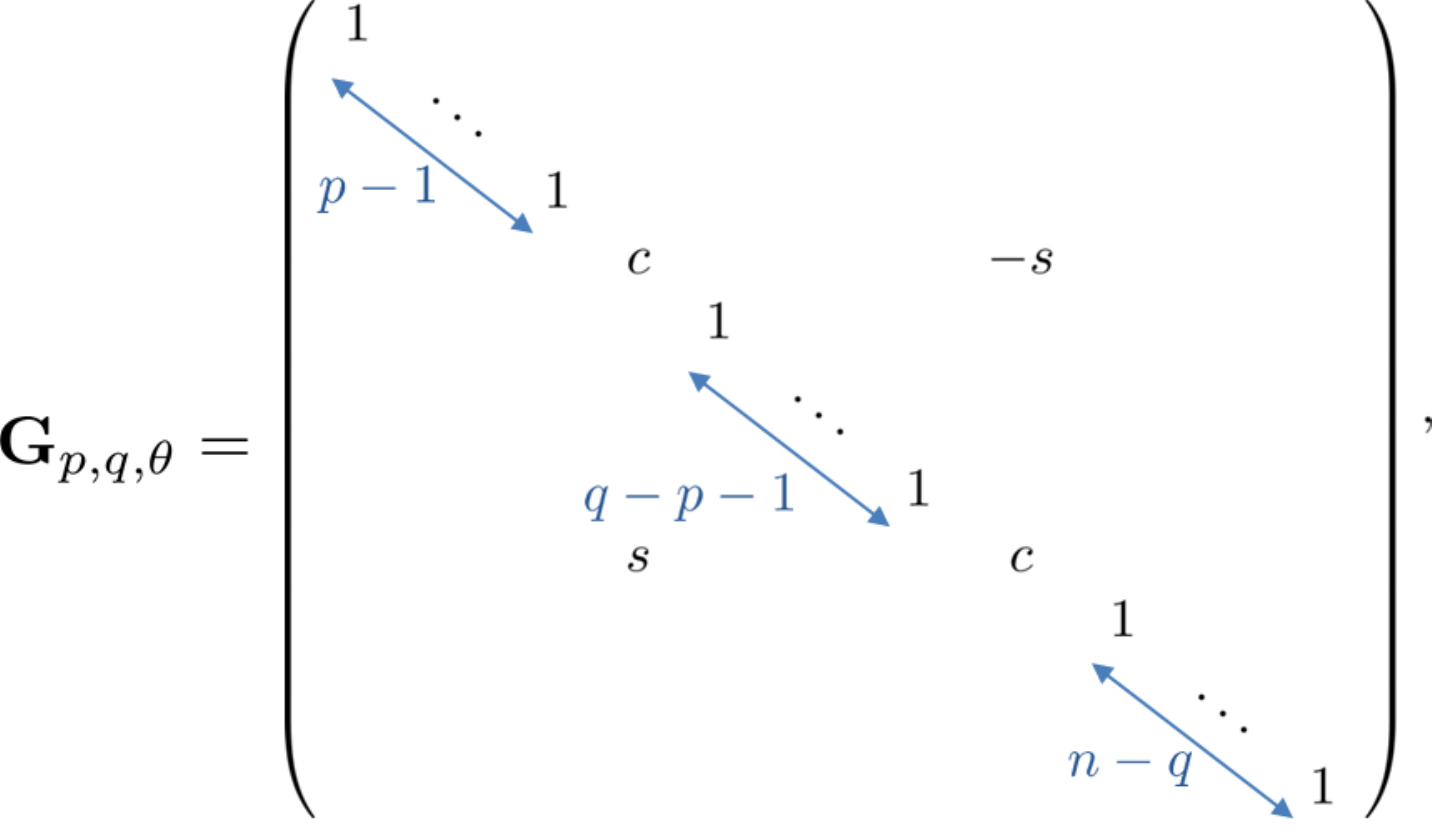}
\end{figure}
\noindent where $c\triangleq\cos(\theta)$ and $s \triangleq \sin(\theta)$. A Givens rotation only depends on three parameters: the two coordinates $p$ and $q$ and the rotation angle $\theta$, hence the notation $\mathbf{G}_{p,q,\theta}$. 
\\

\subsection{Objective}
Our goal is to approximately diagonalize the Laplacian $\mathbf{L}$ with an efficient approximate eigenvector matrix $\hat{\mathbf{U}} = \mathbf{S}_1\dots \mathbf{S}_J$, where the factors $\mathbf{S}_1,\dots,\mathbf{S}_J \in \mathbb{R}^{n \times n}$ are Givens rotations. Using the Frobenius norm to measure the quality of approximation, this objective can be stated as the following optimization problem:
\begin{equation}\tag{DP}
\begin{array}{cl}
\underset{\hat{\boldsymbol{\Lambda}}, \mathbf{S}_1,\dots,\mathbf{S}_J}{\text{minimize}}& \quad \left\Vert \mathbf{L} - \mathbf{S}_1\dots \mathbf{S}_J \hat{\boldsymbol{\Lambda}} \mathbf{S}_J^T\dots \mathbf{S}_1^T \right\Vert_F^2  
\\
& \quad + \sum_{j=1}^J \delta_{\mathcal{S}}(\mathbf{S}_j) + \delta_{\mathcal{D}}(\hat{\boldsymbol{\Lambda}}),
\end{array}
\label{eq:diago}
\end{equation}
where $\mathcal{D}$ is the set of diagonal matrices and $\mathcal{S}$ is the set of Givens rotations.

\subsection{Optimization framework}

{\noindent\textbf{The truncated Jacobi algorithm.}} 
%{\bf RG: since the strategy IS truncated Jacobi, we can't say we propose it}
%\noindent {\bf Greedy strategy.} 
To approximately solve problem~\eqref{eq:diago}, one can rely on a truncated version of the classical Jacobi eigenvalues algorithm \cite{Jacobi1846,Golub2000}. The Jacobi algorithm is an iterative procedure where at each step one seeks the Givens rotation reducing the most the cost function. At the first step, this means setting $\mathbf{S}_1$ and $\hat{\boldsymbol{\Lambda}}$ as follows:
\begin{equation*}
(\mathbf{S}_1,\hat{\boldsymbol{\Lambda}}) \leftarrow \underset{\mathbf{D} \in \mathcal{D}, \mathbf{S} \in \mathcal{S}}{\text{argmin}} \quad \left\Vert \mathbf{L} - \mathbf{S} \mathbf{D}  \mathbf{S}^T \right\Vert_F^2,
\end{equation*}
which can be reformulated, given that the Frobenius norm is invariant under orthogonal transformations, as
\begin{equation*}
(\mathbf{S}_1,\hat{\boldsymbol{\Lambda}}) \leftarrow \underset{\mathbf{D} \in \mathcal{D}, \mathbf{S} \in \mathcal{S}}{\text{argmin}} \quad \left\Vert \mathbf{S}^T\mathbf{L} \mathbf{S}- \mathbf{D}   \right\Vert_F^2.
\end{equation*}
Since $\mathcal{D}$ is the set of all diagonal matrices, the optimal estimated generalized frequencies factor is simply $\hat{\boldsymbol{\Lambda}} = \text{diag}(\mathbf{S}^T\mathbf{L} \mathbf{S})$. This allows to rule out this factor of the problem and to reformulate it as follows:
\begin{equation*}
\mathbf{S}_1 \leftarrow\underset{ \mathbf{S} \in \mathcal{S}}{\text{argmin}} \quad \left\Vert \mathbf{S}^T\mathbf{L} \mathbf{S}   \right\Vert_{\text{offdiag}}^2,
\end{equation*} 
where $\left\Vert \mathbf{A}\right\Vert_{\text{offdiag}}^2$  is the sum of the squared off-diagonal entries of $\mathbf{A}$. Once the factor $\mathbf{S}_1$ is set this way, and introducing the notation $\mathbf{L}_2 \triangleq \mathbf{S}_1^T\mathbf{L} \mathbf{S}_1$, the next step of the strategy is to choose  $\mathbf{S}_2 \leftarrow\underset{ \mathbf{S} \in \mathcal{S}}{\text{argmin}} \quad \left\Vert \mathbf{S}^T\mathbf{L}_2 \mathbf{S}   \right\Vert_{\text{offdiag}}^2$, and so on until the $J$th and last step. The algorithm thus amounts to solve a sequence of $J$ very similar subproblems of the form 
\begin{equation}\tag{SP}
\underset{ \mathbf{S} \in \mathcal{S}}{\text{minimize}} \quad \left\Vert \mathbf{S}^T\mathbf{L}_j \mathbf{S}   \right\Vert_{\text{offdiag}}^2,
\label{eq:subprob}
\end{equation}
with $\mathbf{L}_j \triangleq \mathbf{S}_{j-1}^T \mathbf{L}_{j-1}\mathbf{S}_{j-1}$. This is summarized by the algorithm of {Figure}~\ref{algo:givens_sum}. Compared to the traditional Jacobi eigenvalues algorithm \cite{Jacobi1846,Golub2000}, where new Givens rotations are chosen until a certain accuracy is attained, the main difference is that by prescribing the number $J$ of Givens rotations we can adjust the trade-off between accuracy and computational efficiency of the product $\mathbf{S}_1\dots\mathbf{S}_J$.

\begin{remark}\label{rmk:order}
Note that in order for the approximate FGFT $\hat{\mathbf{U}} = \mathbf{S}_1\dots\mathbf{S}_J$ to make sense in a graph signal processing context, we reorder its columns according to the estimated eigenvalues (line $7$ of the algorithm). This way, the first columns of $\hat{\mathbf{U}} = \mathbf{S}_1\dots\mathbf{S}_J$ ``physically'' correspond to low frequencies, and its last columns to high frequencies.
\end{remark}

\begin{figure}[t]
\centering% for approximate diagonalization
\begin{boxedalgorithmic}[1] 
\REQUIRE{matrix $\mathbf{L}$, target number $J$ of Givens rotations.}\STATE $\mathbf{L}_1 \leftarrow \mathbf{L}$
\FOR{$j=1$ to $J$} 
\STATE  $\mathbf{S}_j \leftarrow \underset{\mathbf{S}\in \mathcal{S}}{\text{argmin}} \left\Vert \mathbf{S}^T\mathbf{L}_j\mathbf{S} \right\Vert_{\text{offdiag}}^2$
\STATE  $\mathbf{L}_{j+1} \leftarrow \mathbf{S}_j^T\mathbf{L}_j\mathbf{S}_j$
\ENDFOR
\STATE  $\hat{\boldsymbol{\Lambda}} \leftarrow \text{diag}(\mathbf{L}_{J+1})$
{\STATE Sort diagonal entries of $\hat{\boldsymbol{\Lambda}}$ in increasing order. Reorder columns of $\mathbf{S}_J$ accordingly.}
\ENSURE sparse orthogonal factors $\mathbf{S}_1,\ldots,\mathbf{S}_J$; diagonal factor $\hat{\boldsymbol{\Lambda}}$.
%\ENSURE The approximate FFT $\hat{\mathbf{U}}_{\text{Givens}} \triangleq \mathbf{S}_1\dots\mathbf{S}_J$ and the aproximate {squared} graph frequencies $\hat{\boldsymbol{\Lambda}}$.
\end{boxedalgorithmic}
\caption{{\bf Truncated Jacobi algorithm}: Approximate diagonalization algorithm with prescribed complexity.}
\label{algo:givens_sum}
\end{figure}

\begin{figure}[tbp]
\raggedright 
\begin{boxedalgorithmic}[1] 
\REQUIRE{matrix $\mathbf{L}_j$.}
\STATE $(p,q) \leftarrow \underset{(r,s)\in [n]^2}{\text{argmax}} |l^j_{rs}|$ 
\STATE $\theta \leftarrow \frac{1}{2}\arctan(\frac{l^j_{qq} - l^j_{pp}}{2l^j_{pq}}) + \frac{\pi}{4}$
\STATE  $\mathbf{S}_j \leftarrow \mathbf{G}_{p,q,\theta} $
\ENSURE matrix $\mathbf{S}_j = \underset{\mathbf{S}\in \mathcal{S}}{\text{argmin}} \left\Vert \mathbf{S}^T\mathbf{L}_j\mathbf{S} \right\Vert_{\text{offdiag}}^2$.
\end{boxedalgorithmic}
\caption{{\bf Resolution of subproblem~\eqref{eq:subprob}}}
\label{algo:subprob}
\end{figure}

{\noindent\textbf{Subproblem resolution.}} 
%\noindent {\bf Subproblem resolution.} 
The algorithm requires to solve $J$ times the optimization subproblem~\eqref{eq:subprob} (at line $3$ of the algorithm of {Figure}~\ref{algo:givens_sum}). The solution of this subproblem is given by the Givens rotation $\mathbf{G}_{p,q,\theta}$, where the indices $p$ and $q$ correspond to the greatest entry of $\mathbf{L}_j$ in absolute value (denoted $|l_{pq}^j|$), and the rotation angle has the expression $\theta = \frac{1}{2}\arctan(\frac{l^j_{qq} - l^j_{pp}}{2l^j_{pq}}) + (2k+1)\frac{\pi}{4}$, $k \in \mathbb{Z}$. We then have $\left\Vert \mathbf{L}_{j+1} \right\Vert_{\text{offdiag}}^2 = \left\Vert \mathbf{L}_j \right\Vert_{\text{offdiag}}^2 - 2(l^j_{pq})^2 $. For a proof as well as a review of the different implementations and variants of the Jacobi algorithm, see \cite[pages 426-435]{Golub2012}.  Figure~\ref{algo:subprob} details the algorithm to solve subproblem~\eqref{eq:subprob}. 
\\

\noindent{\textbf{Refinements.}} This truncated Jacobi algorithm requires $\mathcal{O}(n^2 + nJ)$ operations to compute. We proposed in~\cite{magoarou_approximate_2017} a parallel framework that reduces this complexity to $\mathcal{O}(nJ\log{n})$. Nevertheless, in this paper, we will not concentrate on the cost of \emph{obtaining} the sparse matrix factorization; we will rather concentrate on the tradeoff between the cost of \emph{applying} $\hat{\mathbf{U}}$ to a graph signal, which, in both cases (truncated or parallel truncated Jacobi), is $\mathcal{O}(J)$,  and the quality of the approximation of $\mathbf{U}$ by $\hat{\mathbf{U}}$. The codes of this truncated Jacobi algorithm are freely available at \texttt{https://faust.inria.fr/}. 

\section{Experiments}
In the following, let us denote by $\mathbf{u}_k$ the $k$-th eigenvector of $\mathbf{L}$; i.e., $\mathbf{U}=(\mathbf{u}_1|\ldots|\mathbf{u}_n)$. Let us also denote by $\hat{\mathbf{u}}_k^{J}$ the $k$-th approximate eigenvector of $\mathbf{L}$; i.e., $\hat{\mathbf{U}}= \mathbf{S}_1\dots\mathbf{S}_J =(\hat{\mathbf{u}}_1^{J}|\ldots|\hat{\mathbf{u}}_n^{J})$. 
Given this approximate fast graph Fourier matrix $\hat{\mathbf{U}}$, a natural question is to investigate how well does this matrix approximate the real graph Fourier matrix $\mathbf{U}$ as a function of $J$, the number of Givens rotations\footnote{Note that the orientation of $\mathbf{u}_i^J$ is arbitrary: both $\hat{\mathbf{u}}_i^J$ and $-\hat{\mathbf{u}}_i^J$ are equivalent from a diagonalisation point of view. For a fair comparison with $\mathbf{u}_i$, the orientation of $\hat{\mathbf{u}}_i^J$ is decided by settling that $\mathbf{u}_i^\top \hat{\mathbf{u}}_i^J$ should be positive.}. To conduct this analysis, we propose the following experiments. 

\subsection{Types of graphs} 
We consider three types of random graphs:
\begin{itemize}
 \item {\bf Erdos-Renyi graphs}. We set the number of nodes $n$ to $128$. All pairs of nodes are connected with probability $p$. The average degree of such a model is $c=p(n-1)$. 
 \item {\bf Stochastic Block Model (SBM)}, a random community-structured graph model. We specifically look at graphs 
with $m$ communities of same size $n/m$. In the SBM, the probability of connection between any two nodes $i$ and $j$ is $q_1$ if they are in the same community, and $q_2$ otherwise. One can show that the average degree reads  $c=q_1\left(\frac{n}{k}-1\right)+q_2\left(n-\frac{n}{k}\right)$. 
Thus, instead of providing the probabilities $(q_1,q_2)$, one may characterize a SBM by considering $(\epsilon=\frac{q_2}{q_1},c)$. 
The larger $\epsilon$, the fuzzier the community structure. In fact, 
authors in~\cite{decelle_asymptotic_2011} show that above the critical value $\epsilon_c=(c-\sqrt{c})/(c+\sqrt{c}(m-1))$, community structure becomes  undetectable in the large $n$ limit. In the following, we set $n=128$ and $m=8$. %The average degree $c$ and the parameter $\epsilon$ will be varying.  
 \item {\bf Sensor graphs.} $n=128$ sensors (points in two dimensions) are randomly and uniformly generated on the unit square. Two sensors $i$ and $j$ are connected if their interdistance is inferior to a given threshold $\tau$. The obtained graph is binary. The average degree $c$ of such a model is not explicit. 
\end{itemize}
%We also consider two real-world graphs:
%\begin{itemize}
% \item the Minnesota road network. 
% \item the Bunny graph. 
%\end{itemize}
For the sensor graph, we set $\tau=0.161$, which yields graphs with average degree $c=10$. In experiments with Erdos-Renyi and SBM graphs, we enforce the same average degree by setting $p=c/(n-1)=10/127$ for the Erdos-Renyi model, and $c=10$ for the SBM model (besides, $\epsilon$ is set to $\epsilon_c/10$).

\begin{figure*}[htbp]
 \centering
a)\includegraphics[width=0.59\columnwidth]{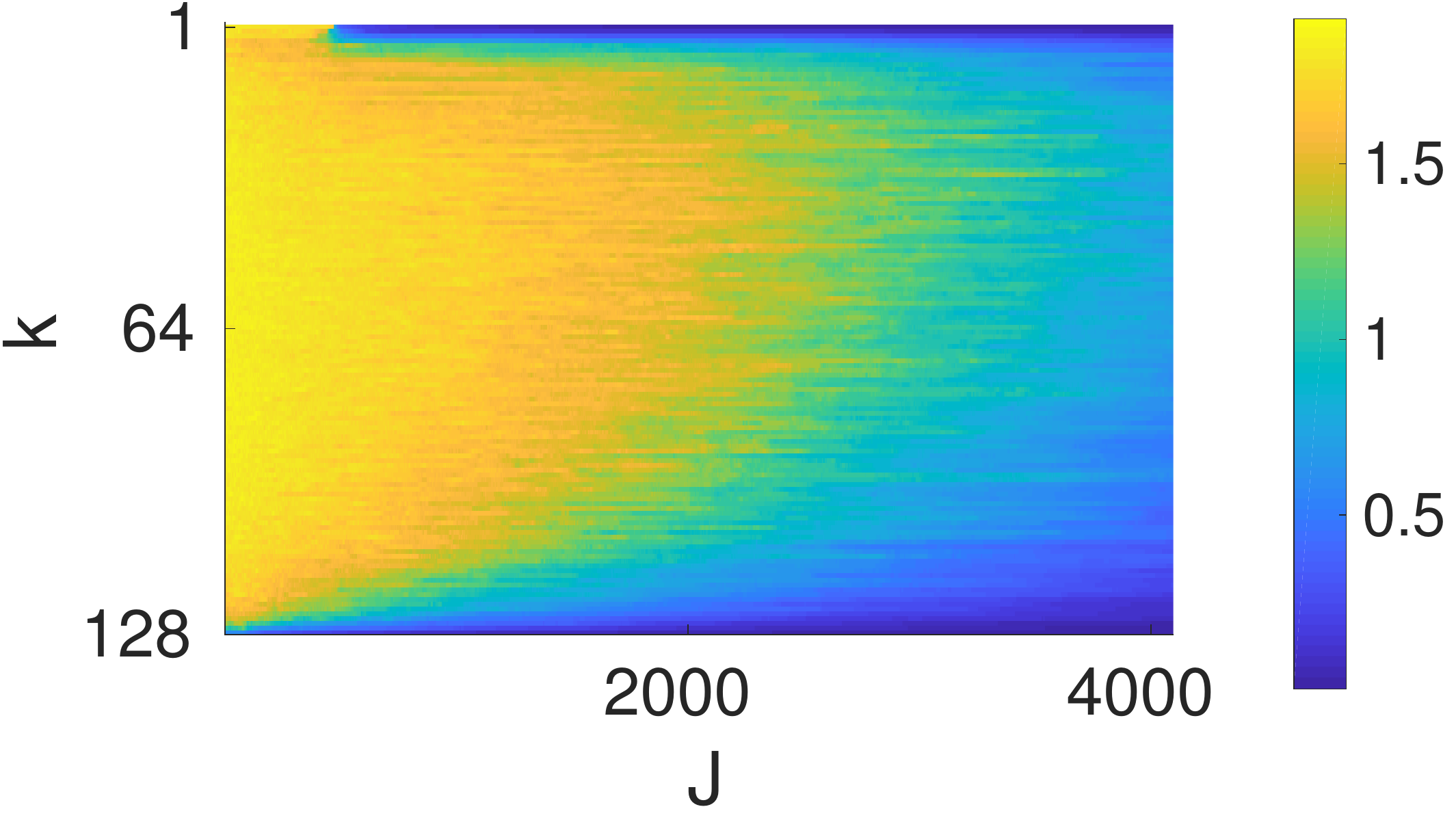} \hfill
b)\includegraphics[width=0.59\columnwidth]{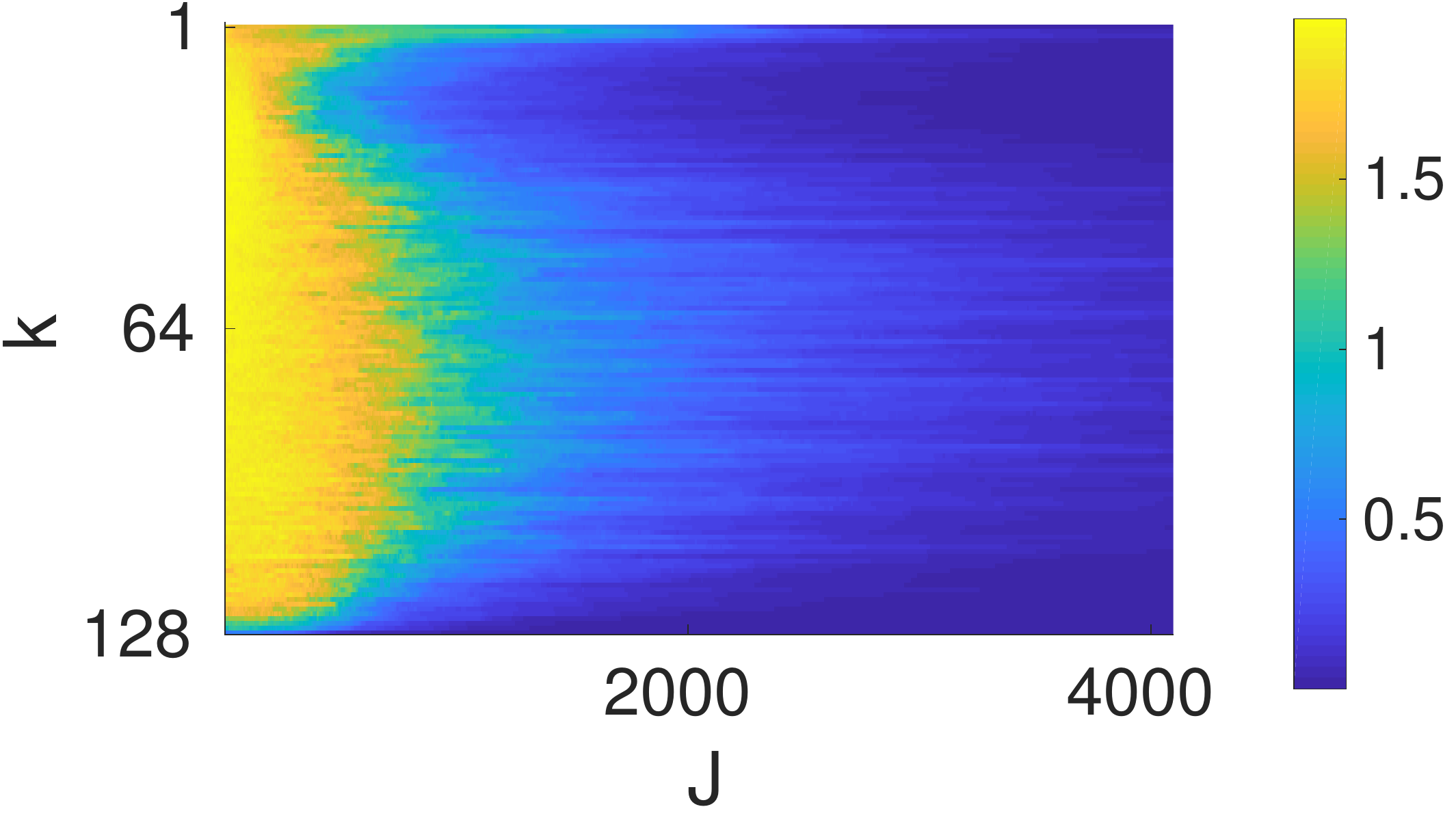}  \hfill
c)\includegraphics[width=0.59\columnwidth]{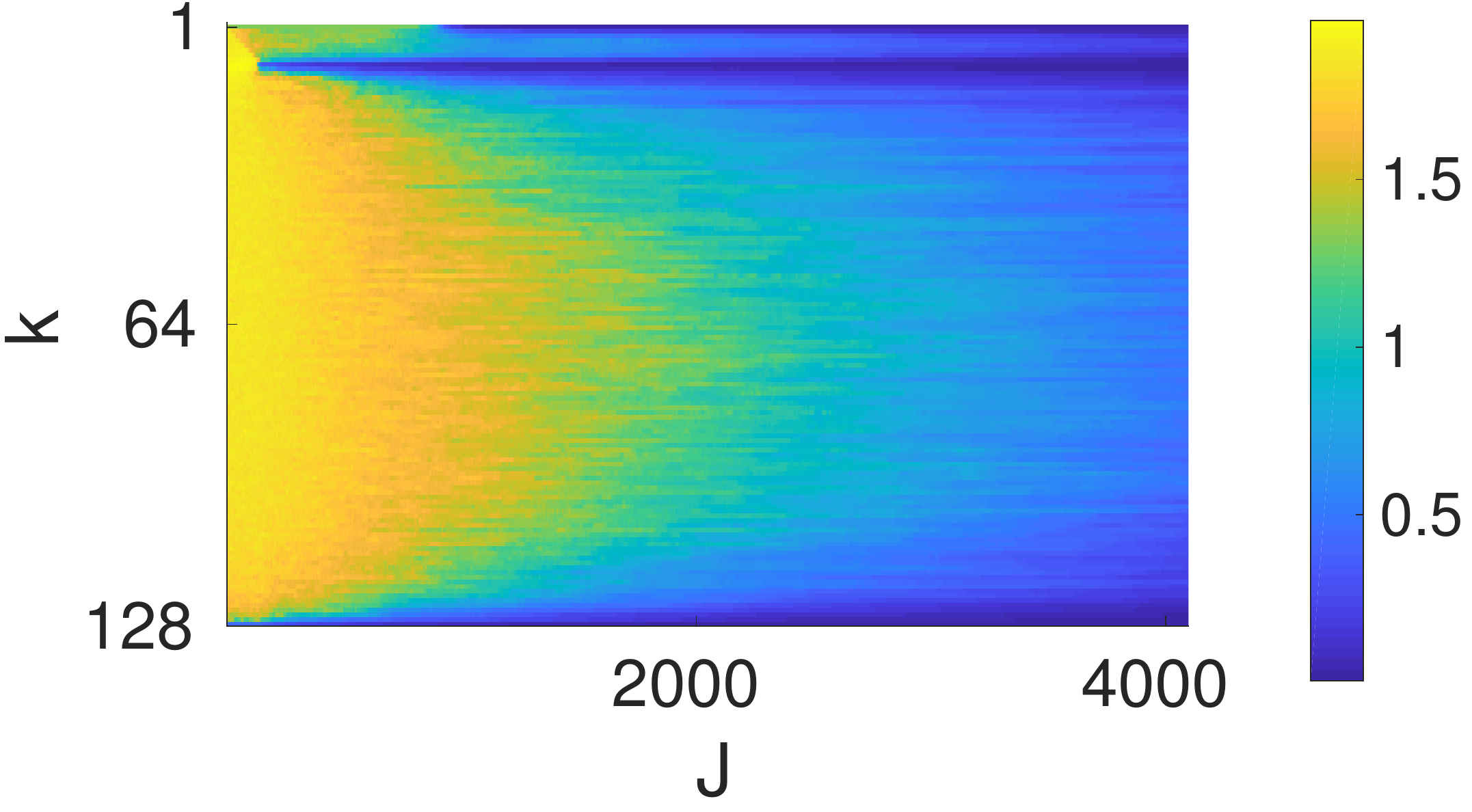} %\\
%~~~~~~~~~~~~~a) ~~~~~~~~~~~~~~~~~~~~~~~~~~~~~~~~~~ b) ~~~~~~~~~~~~~~~~~~~~~~~~~~~~~~~~~ c) ~~~~~~~~~~~~~~~~~~~~~~~~~~~~~~~d) ~~~~~~~~~~~~~~~~~~~~~~
\caption{Median of the error $\text{err}_{1}(k,J)^2$ over 100 random draws: a)~Erdos-Renyi graph, b)~sensor graph, c)~SBM graph. }
\label{fig:experiments_err_c_matrix}
\end{figure*}

%First of all, we decide to fix the number of Givens rotations to $J=n\log{n}$, in order to mimick the classical FFT. 
\subsection{Analysis of the approximation errors}
In~\cite{magoarou_approximate_2017}, we concentrated on global approximation error measures such as:
\begin{equation}\label{eq:GlobalError}
%\text{err}_1(J) \triangleq 
\frac{\Vert \mathbf{U} - \hat{\mathbf{U}} \Vert_F}{\left\Vert \mathbf{U} \right\Vert_F}.
\end{equation}
This measure does not differentiate \emph{where} in the spectrum are the errors. The goal of this paper is to refine this  analysis. 
First of all, note that:
 \begin{align*}
  %\text{err}_1(J)^2 &= 
  \frac{\Vert \mathbf{U} - \hat{\mathbf{U}} \Vert_F^2}{\left\Vert \mathbf{U} \right\Vert_F^2}
  &=\frac{\Vert \mathbf{Id} - \mathbf{U}^\top\hat{\mathbf{U}} \Vert_F^2}{\left\Vert \mathbf{U} \right\Vert_F^2}\\
  &=\frac{1}{\left\Vert \mathbf{U} \right\Vert_F^2}\sum_{k=1}^n\sum_{j=1}^n (\delta_{jk} - \mathbf{u}_j^\top\hat{\mathbf{u}}_k^{J})^2\\
  &= \frac{1}{n}\sum_{k=1}^n \text{err}_{1}(k,J)^2,
 \end{align*}
 where we consider the fine-grain error measure
 \begin{equation}\label{eq:FineGrainError1}
 \text{err}_{1}(k,J)^2 := \left\Vert\bm{\delta}_k-\mathbf{U}^\top\hat{\mathbf{u}}_k^{J}\right\Vert_2^2
 \end{equation}
that enables a finer analysis of the approximation error~\eqref{eq:GlobalError}.

 \begin{figure*}
 \centering
a)\includegraphics[width=0.59\columnwidth]{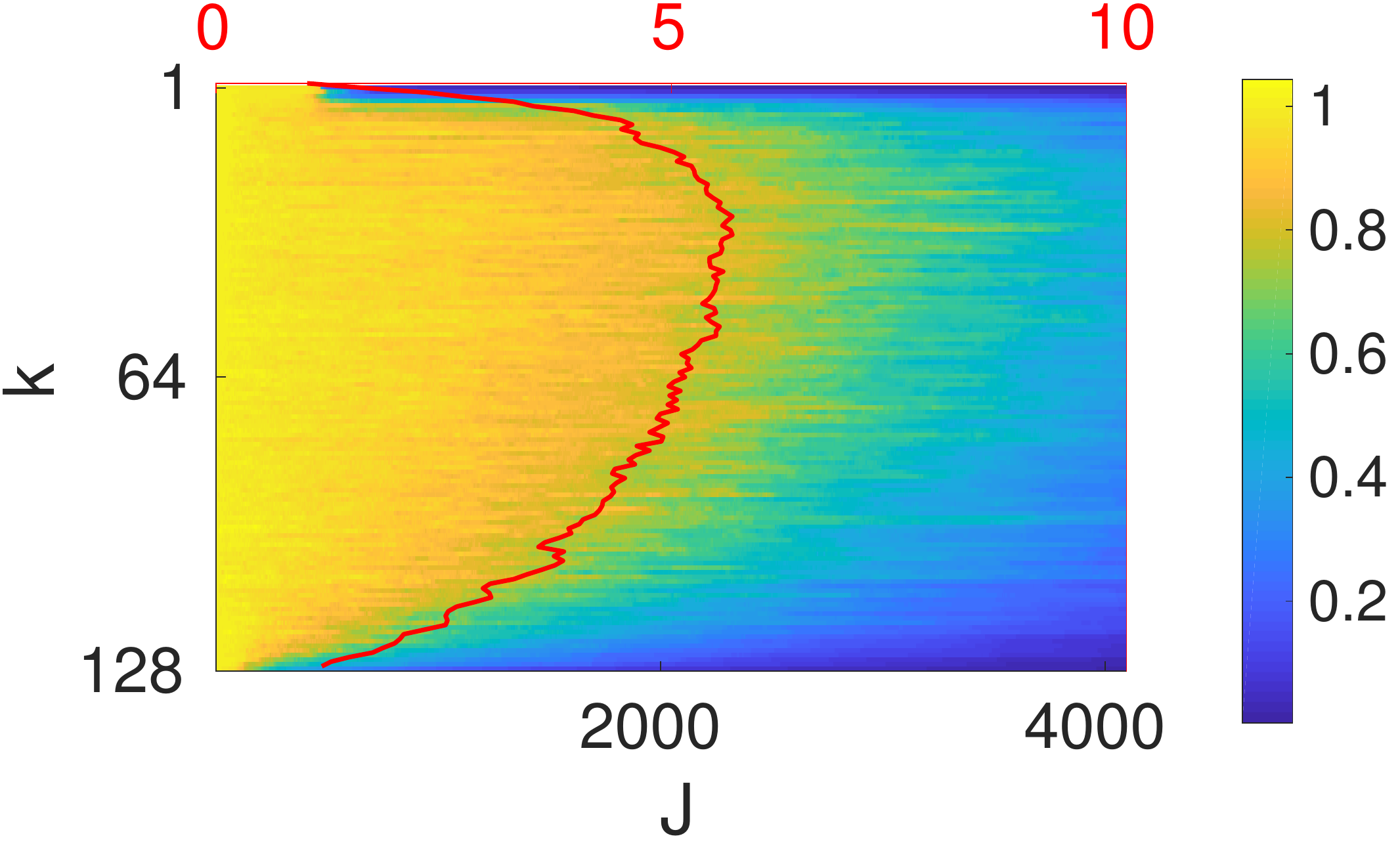} \hfill
b)\includegraphics[width=0.59\columnwidth]{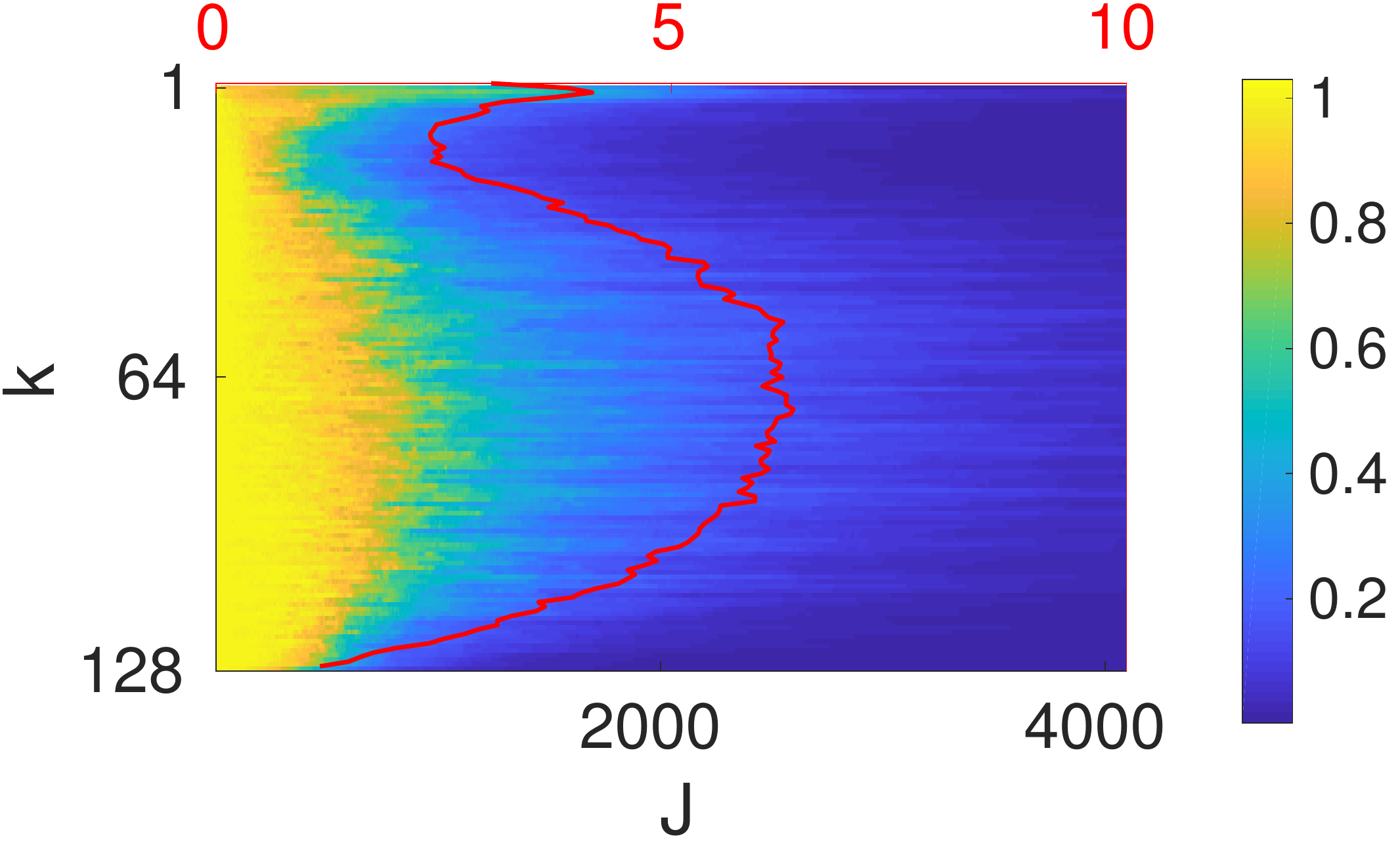}  \hfill
c)\includegraphics[width=0.59\columnwidth]{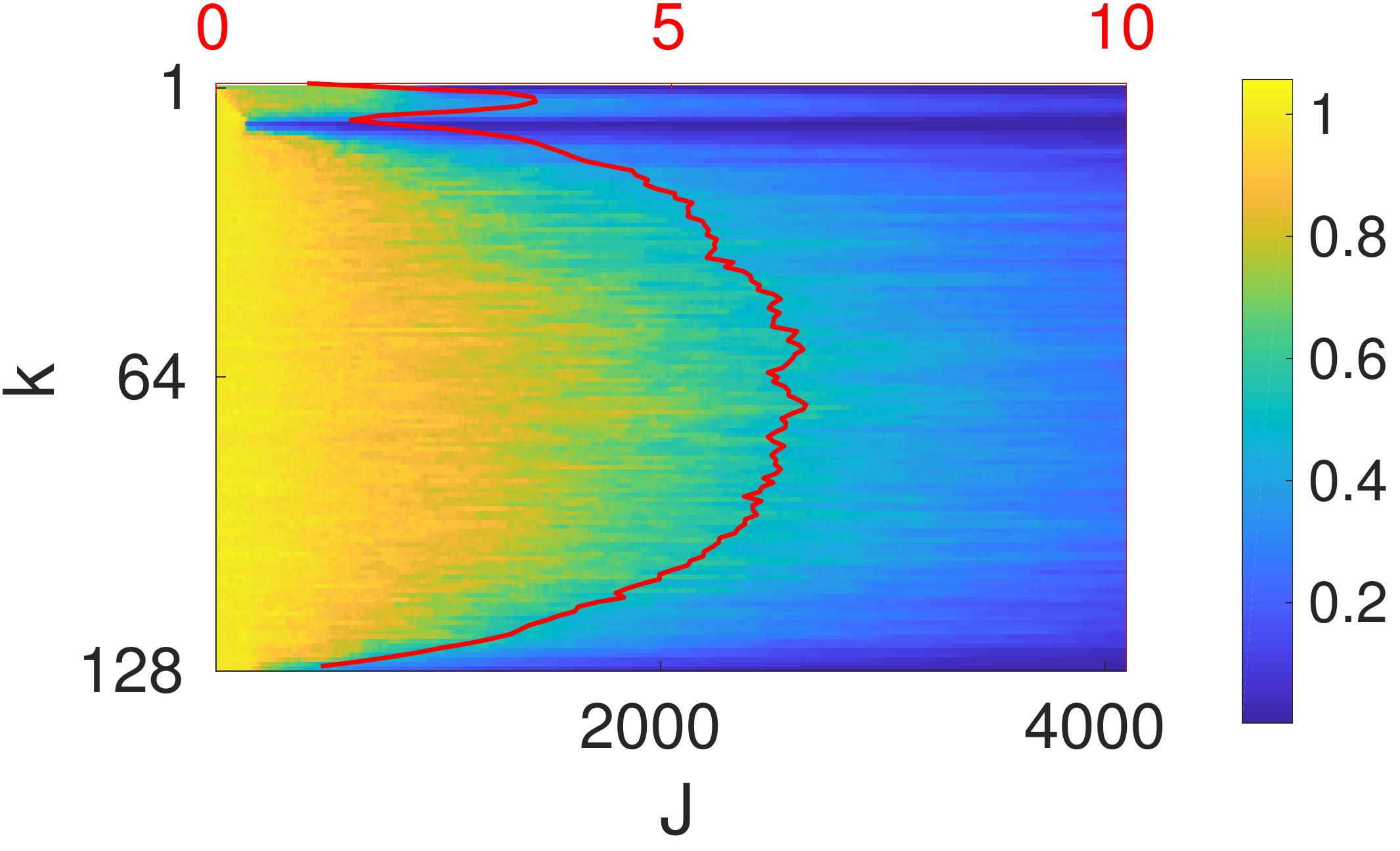} %\\
%~~~~~~~~~~~~~a) ~~~~~~~~~~~~~~~~~~~~~~~~~~~~~~~~~~ b) ~~~~~~~~~~~~~~~~~~~~~~~~~~~~~~~~~ c) ~~~~~~~~~~~~~~~~~~~~~~~~~~~~~~~d) ~~~~~~~~~~~~~~~~~~~~~~
\caption{
Median of the (degree-corrected) normalized error $\widetilde{\text{err}}_{1}(k,J)^2$ over 100 random draws: a)~Erdos-Renyi graph, b)~sensor graph, c)~SBM graph. The red line represents a measure of eigenvalue density as a function of $k$ (see Eq.~\eqref{eq:eig_density}).}
\label{fig:experiments_err_c_matrix_norm}
\end{figure*}

In Figure~\ref{fig:experiments_err_c_matrix}, we show $\text{err}_{1}(k,J)^2$ as a function of $k$ and $J$ for the three random models detailed above. We readily observe very different behaviors. Firstly, for small values of $J$, the approximation error is large and somewhat uniform over the whole spectrum, except for high values of $k$ where the error drops. Moreover, as $J$ increases, the approximation error decreases, but not uniformly: the error decreases faster for some values of $k$. In the Erdos-Renyi example, the approximation error decreases faster on the boundaries (small or large values of $k$) than in the middle. In the random sensor graph example, the approximation error decreases faster around $k=10$ and for large values of $k$. For the SBM graph, the error decreases very fast around $k=8$ (notice that this coincides with $k=m$, the number of communities), and faster on the boundaries than in the middle. In the following, we give tentative explanations accounting for these differences. 
\\

\noindent\textbf{Normalized error measure.} For $J=0$ we have $\hat{\mathbf{U}}=\mathbf{Id}$, up to column permutation (cf Remark~\ref{rmk:order}) according to the diagonal of $\hat{\mathbf{U}}^\top\mathbf{L}\hat{\mathbf{U}}=\mathbf{L}$, i.e. with respect to $\mathbf{D}$, the degree matrix. Taking into account this ordering yields
\begin{align}\label{eq:DegreeOrder}
 \hat{\mathbf{U}}=\left(\bm{\delta}_{\sigma(1)}|\bm{\delta}_{\sigma(2)}|\ldots|\bm{\delta}_{\sigma(n)}\right),
\end{align}
with $\sigma$ a permutation such that 
\(
d_{\sigma(1),\sigma(1)}\leq d_{\sigma(2),\sigma(2)}\leq \ldots \leq d_{\sigma(n),\sigma(n)}.
\)
In particular $\sigma(n)$ indexes the node with largest degree, and one has:
 $$\text{err}_{1}(k,0)^2 = \left\Vert\bm{\delta}_k-\mathbf{U}^\top\bm{\delta}_{\sigma(k)}\right\Vert_2^2.$$

\begin{figure}
 \centering
\includegraphics[width=0.61\columnwidth]{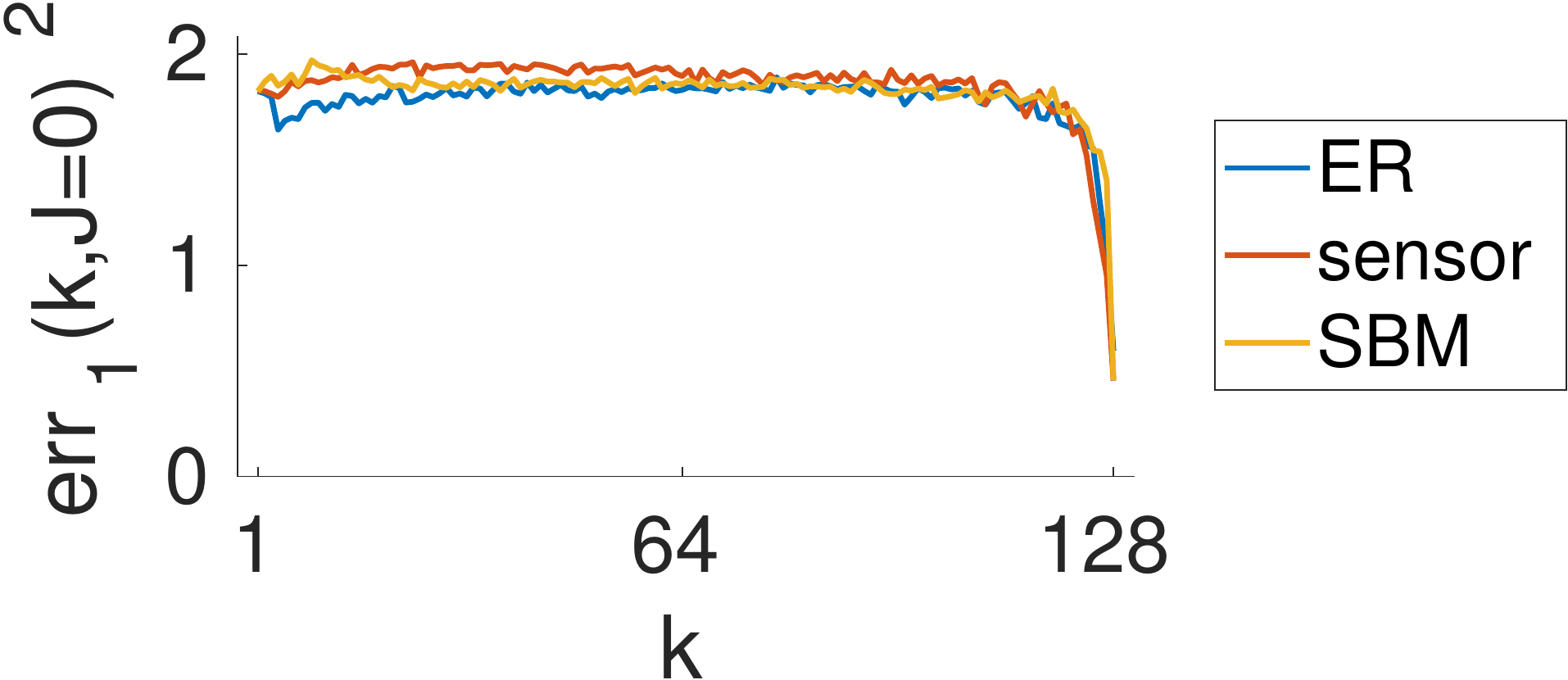} \\\vspace{0.3cm}
\includegraphics[width=0.61\columnwidth]{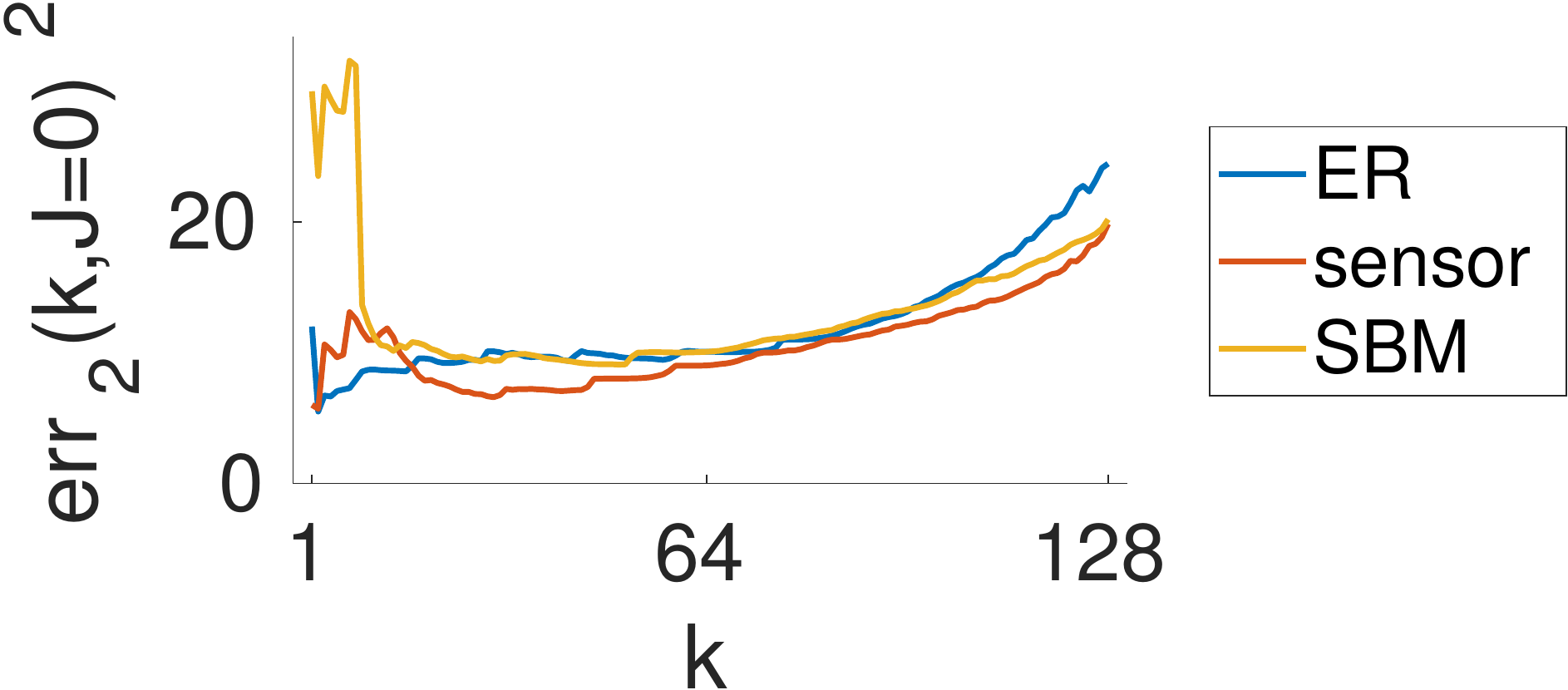} 
%~~~~~~~~~~~~a) ~~~~~~~~~~~~~~~~~~~~~~~~~~~~~~~~~~ b) ~~~~~~~~~~~~~~~~~~~~~~~~~~~~~~~~~ c) ~~~~~~~~~~~~~~~~~~~~~~~~~~~~~~~d) ~~~~~~~~~~~~~~~~~~~~~~
\caption{Median error $\text{err}_{1}(k,0)^2$ (top) [resp. $\text{err}_{2}(k,0)^2$ (bottom)] over 100 draws.}
\label{fig:experiments_err_0}
\end{figure}

As shown on Figure~\ref{fig:experiments_err_0}-a, this error is not evenly distributed along the spectrum:  higher values of $k$ show a comparatively lower approximation error. As $\mathbf{U}^\top\bm{\delta}_{\sigma(k)}=(u_1(\sigma(k)),\ldots,u_n(\sigma(k)))^\top$, this non-uniformity shows in fact that high degree nodes tend to contribute more to high-frequency Fourier modes than other nodes. In other words, high frequency Fourier modes tend to localize more on high degree nodes, at least on the three random graph models considered here. 

To validate that this accounts for the non-uniformity observed for large $k$ and small $J$ on Figure~\ref{fig:experiments_err_c_matrix}, we define a  normalized  error measure as:
\begin{align}\label{eq:DegreeCorrectedError}
 \widetilde{\text{err}}_{1}(k,J)^2 = \frac{\text{err}_{1}(k,J)^2}{\text{err}_{1}(k,0)^2}. 
\end{align}
As shown on Figure~\ref{fig:experiments_err_c_matrix_norm}, the non-uniformity at low values of $J$ is corrected, but the heterogeneity over the spectrum of the decrease of the error with respect to $k$ is still unaccounted for.

 \begin{figure*}
 \centering
a)\includegraphics[width=0.59\columnwidth]{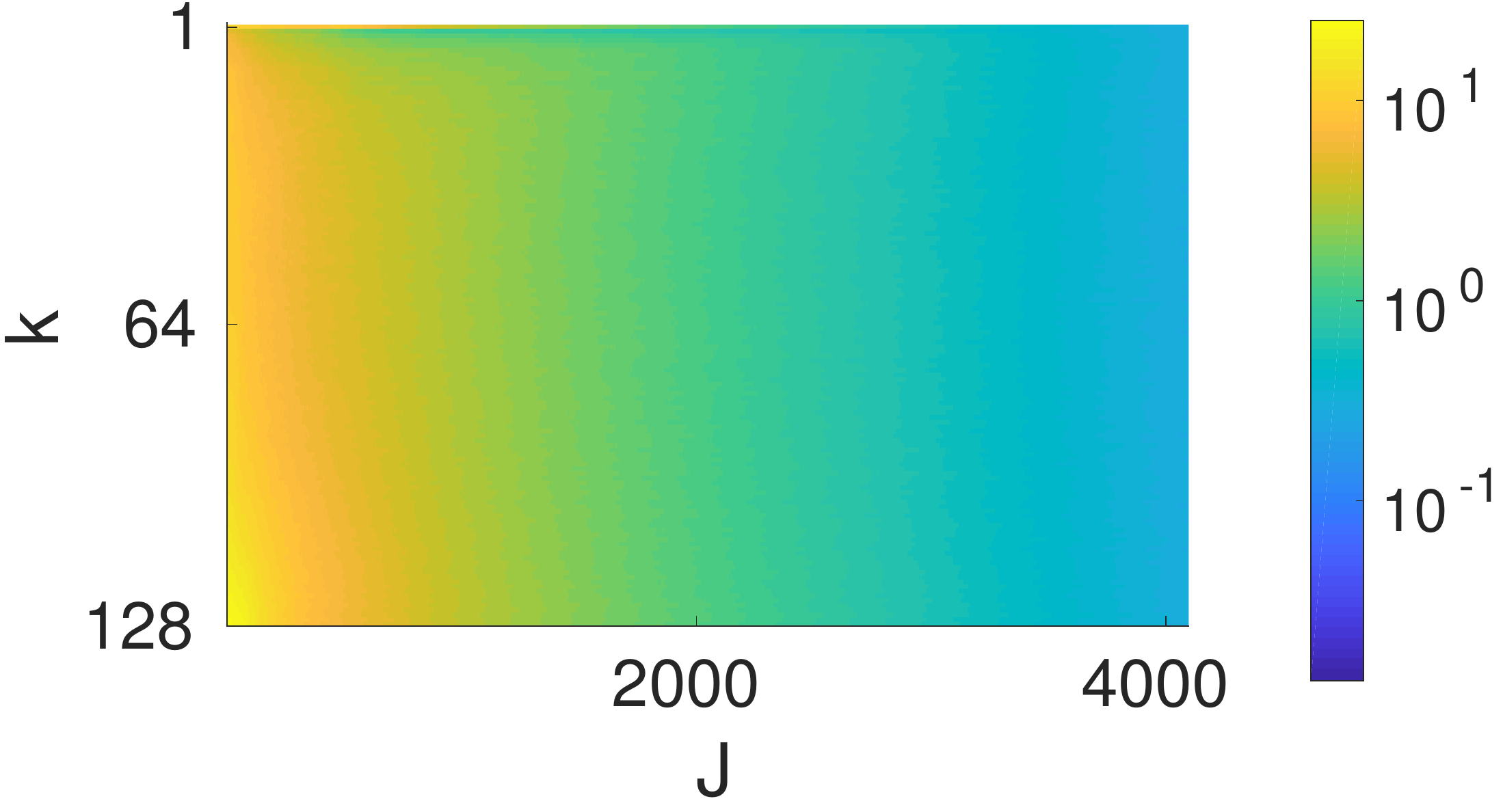} \hfill
b)\includegraphics[width=0.59\columnwidth]{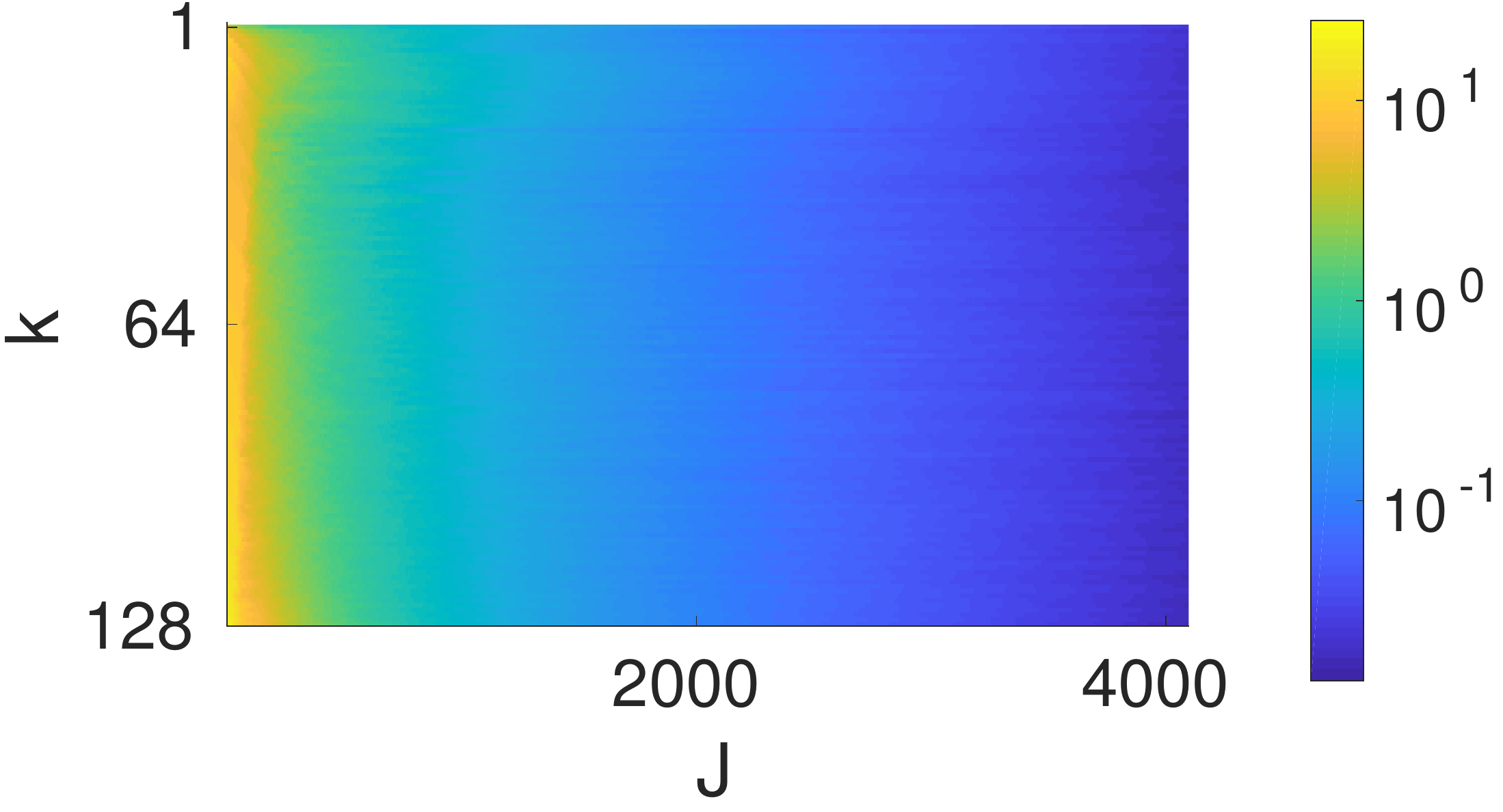}  \hfill
c)\includegraphics[width=0.59\columnwidth]{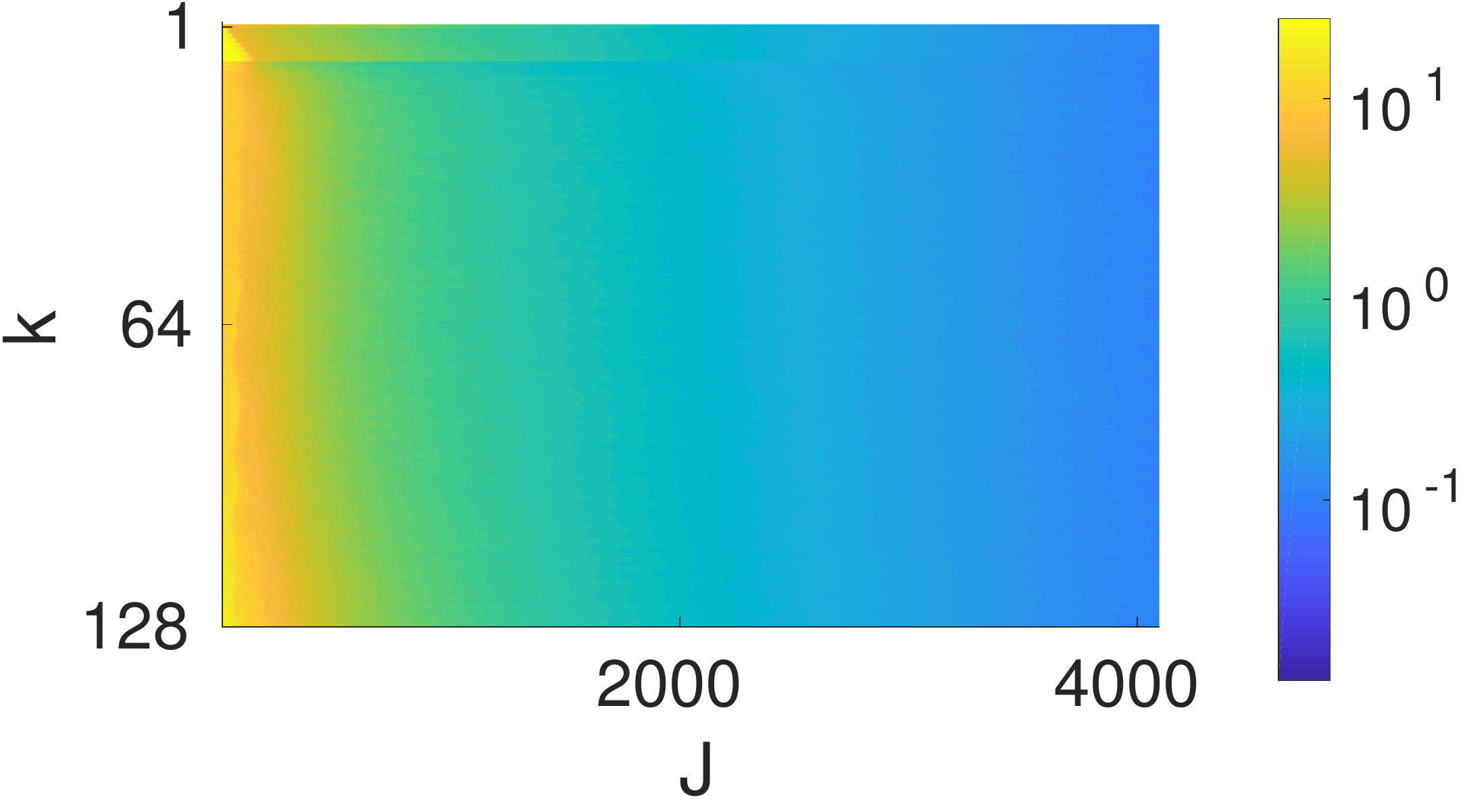} \\\vspace{0.3cm}

d)\includegraphics[width=0.59\columnwidth]{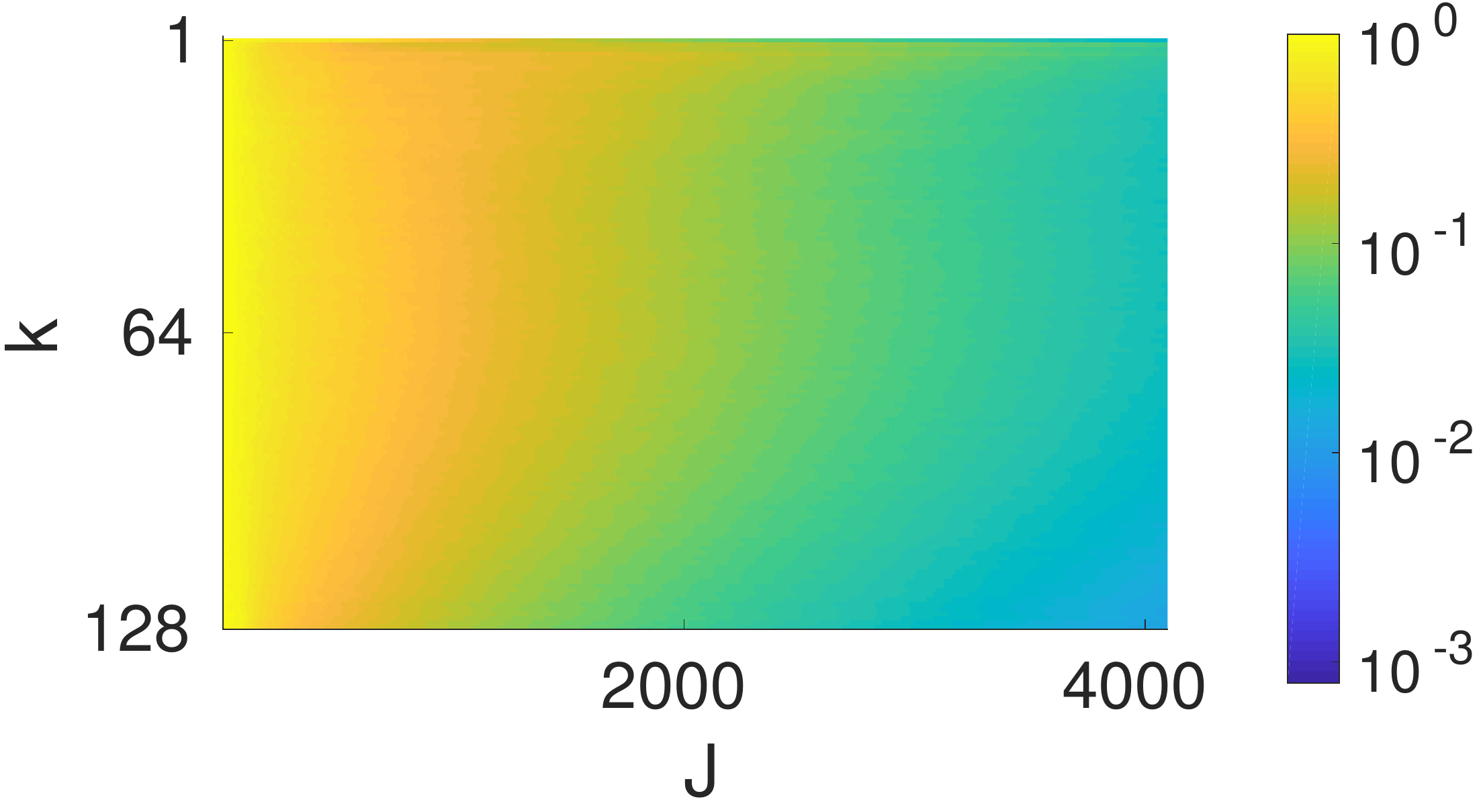} \hfill
e)\includegraphics[width=0.59\columnwidth]{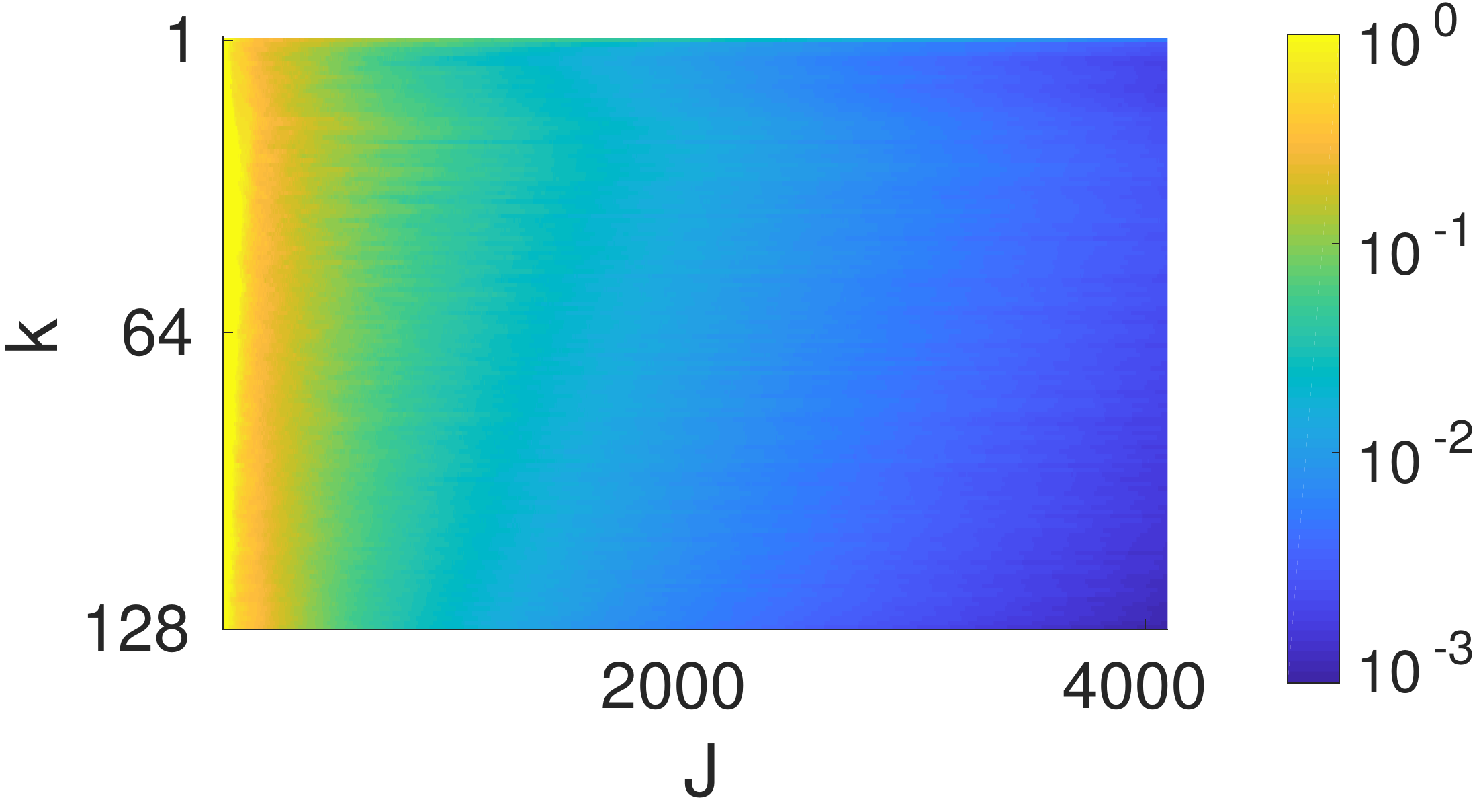}  \hfill
f)\includegraphics[width=0.59\columnwidth]{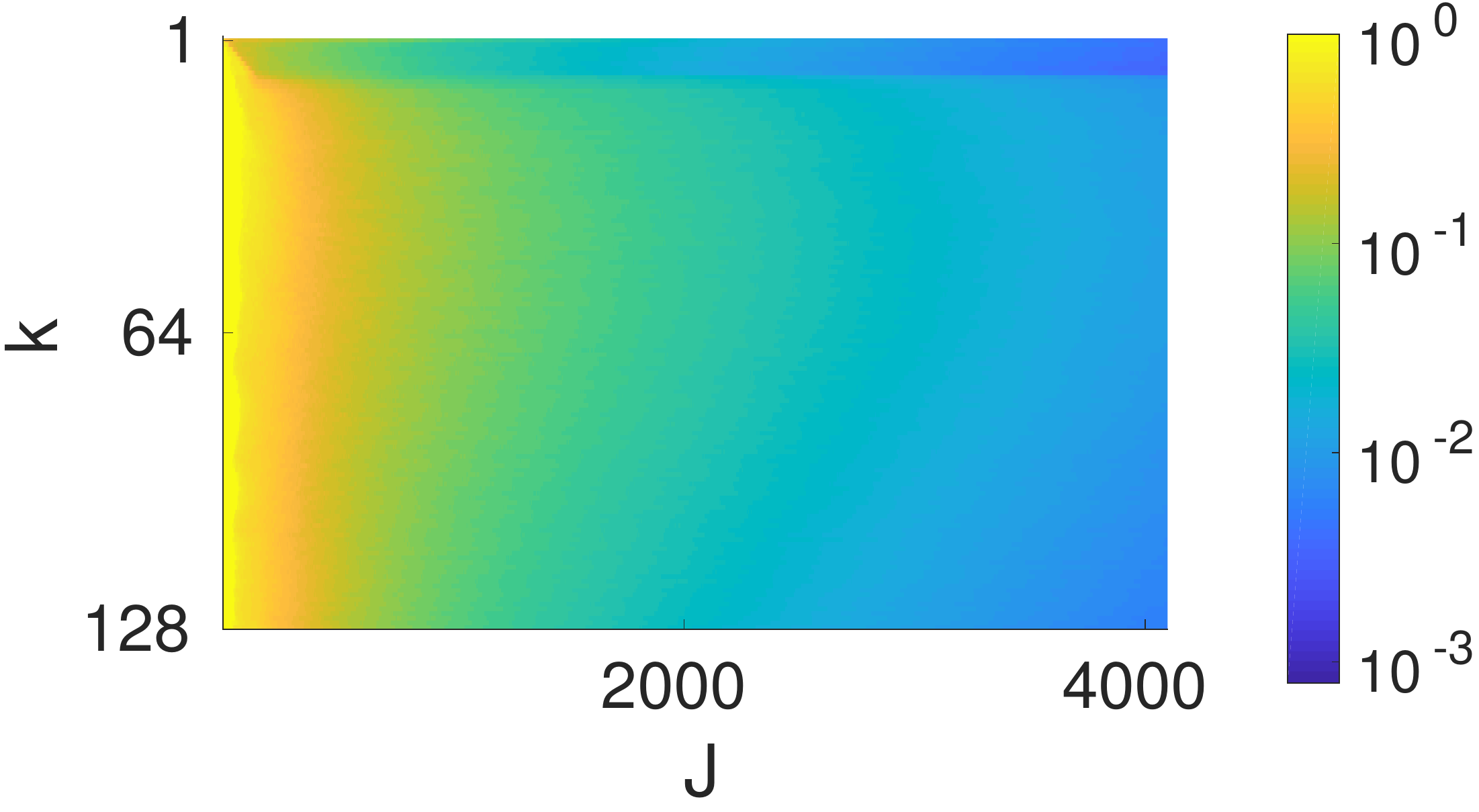} \\
%~~~~~~~~~~~~~a) ~~~~~~~~~~~~~~~~~~~~~~~~~~~~~~~~~~ b) ~~~~~~~~~~~~~~~~~~~~~~~~~~~~~~~~~ c) ~~~~~~~~~~~~~~~~~~~~~~~~~~~~~~~d) ~~~~~~~~~~~~~~~~~~~~~~
\caption{The first (resp. second) line of figures represents $\text{err}_{2}(k,J)^2$ (resp. $\widetilde{\text{err}}_{2}(k,J)^2$). Note that the color scale is logarithmic. Results shown here are the median behaviour over 100 realisations of the underlying random model: a+d)~Erdos-Renyi graph, b+e)~sensor graph, c+f)~SBM graph.}
\label{fig:experiments_err_c_2_matrix_norm}
\end{figure*}

{\bf  Eigenvalue density.}
The remaining heterogeneity can be conjectured to be due to the inhomogeneity of the spectrum density: in the extreme case of eigenvectors associated to an eigenvalue of multiplicity greater than one, the error cannot be small. To elaborate, define an eigenvalue density measure:
 \begin{align}
 \label{eq:eig_density}
  f(k) = \#\{i\text{ s.t. } \lambda_i\in[\lambda_k-\Delta, \lambda_k+\Delta].
 \end{align}
 If $\lambda_k$ is isolated and there are no other eigenvalues in its vicinity (defined by $\pm \Delta$), then $f(k)=1$. Whereas if $\lambda_k$ is in a high-density part of the spectrum, then $f(k)$ is large. Superimposed on Figure~\ref{fig:experiments_err_c_matrix_norm} representing $\widetilde{\text{err}}_1(k,J)$, we plot $f(k)$ (computed with $\Delta = 0.25$) as a function of $k$. We indeed observe an interesting correlation between the approximation error for $k$ and the density of the spectrum around $\lambda_k$: the higher the eigenvalue density, the larger the error.

% \noindent{\textbf{A second error measure.}}
To validate this correlation, we define an error measure that probes how well $\hat{\mathbf{u}}_k^{J}$ approximates $\mathbf{u}_k$ as an eigenvector of $\mathbf{L}$:
\begin{equation}
  \text{err}_2(k,J)^2 = \norm{\Lap\hat{\fou}_k^{J}-\lambda_k\hat{\fou}_k^{J}}_2^2,
\end{equation}
where $\lambda_k$ is the exact $k$-th eigenvalue of $\mathbf{L}$. This measure may be rewritten as:
\begin{align*}
 \text{err}_2(k,J)^2 
 &= \norm{\Fou\Lambda\Fou^\top\hat{\fou}_k^{J}-\lambda_k\hat{\fou}_k^{J}}_2^2
 =\norm{\Lambda\Fou^\top\hat{\fou}_k^{J}-\lambda_k\Fou^\top\hat{\fou}_k^{J}}_2^2\\
 &=\sum_{j=1}^n (\lambda_j-\lambda_k)^2 (\fou_j^\adjoint\hat{\fou}_k^{J})^2.
\end{align*}
To illustrate the difference between the two errors, say  $\hat{\fou}_k^{J}$ is a weighted sum of two true eigenvectors $\hat{\fou}_k^{J} = \alpha\fou_k+\beta\fou_{k+1}$, with $\alpha^2+\beta^2=1$. The first error measures $\text{err}_1^2 = 2(1-\alpha)$ and does not depend on the distance between $\lambda_{k+1}$ and $\lambda_k$, whereas the second error measures $\text{err}_2^2 = (\lambda_{k+1}-\lambda_k)^2(1-\alpha^{2})$, meaning that if $\hat{\fou}_k^{J}$ is not exactly $\fou_k$ but close to it from a ``spectral point-of-view'', then the error is small. Once again, one may normalize this error to correct for degree heterogeneities, using the notations of~\eqref{eq:DegreeOrder} to write:
\begin{align*}
\text{err}_{2}(k,0)^2 = \left\Vert\mathbf{L}\bm{\delta}_{\sigma(k)}-\lambda_k\bm{\delta}_{\sigma(k)}\right\Vert_2^2.
\end{align*}
%This error is a trade-off between the density of the vector $\mathbf{L}\bm{\delta}_{\sigma(k)}$ (the larger $k$, the higher the degree of ${\sigma(k)}$, the denser its associated column in the Laplacian) and the value of $\lambda_k$ (the larger $k$, the larger $\lambda_k$). 

Figure~\ref{fig:experiments_err_0}-b shows how $\text{err}_2(k,0)$ varies with $k$.  We observe how it increases at extreme values of $k$ (very small or very large). Interestingly, for the SBM, the large error at small $k$ is due to the community structure that imposes the first $m=8$ eigenvalues (because there are 8 communities in the model) to be very small. We checked (experiments not displayed here) that decreasing the strength of the community structure (by increasing $\epsilon$ to values closer to $\epsilon_c$) increases the first eigenvalues, and decreases the error at $J=0$. 

 Figure~\ref{fig:experiments_err_c_2_matrix_norm} displays $\text{err}_2(k,J)$ and its normalized version $\widetilde{\text{err}}_2(k,J)$. For visualisation purposes, the color scale is logarithmic. This new measure drastically decreases the inhomogeneity of the error over the spectrum. 

%{\bf Ici on aimerait desormais voir si le $J$ qui permet d'atteindre un $\tilde{\text{err}}_{2}(\star,J) \leq \epsilon$ depend logarithmiquement de la taille du graphe ?}

\section{Conclusion}

In this paper, the approximation capabilities of FGFTs were empirically studied. In order to do so, the FGFT of graphs pertaining to various families were computed and their approximation errors analyzed with respect to the number of used Givens rotations. Two main conclusions can be drawn. First, high frequency modes tend to be localized around high degree nodes, and are thus well approximated with few Givens rotations (that are very sparse matrices). Second, Fourier modes corresponding to dense regions of the spectrum are more difficult to disentangle, and thus need more Givens rotations to be well approximated. Stated differently, if one wishes to exactly recover the ordering of the eigenmodes, FGFTs perform poorly in the dense parts of the spectrum. On the other hand, if one only whishes to approximate eigenspaces in the vicinity of a given eigenvalue, the FGFT is an interesting tool to consider and leads to fast implementations. 
%
%{\bf *** Nico a ecrit: la TF sur graphe avec l'algo de Jacobi contient principalement des erreurs dans des endroits où la densité spectrale est importante (en fait, c'est au niveau de l'ordonnencement qu'il y a des erreurs, parce que les valeurs propres sont trop proches).  Donc si on veut un ordre strict de la TF d'un signal, c'est pas terrible. Par contre si on veut qu'au voisinage de n'importe quelle valeur propre, les espaces propres soient assez bien approchés, alors c'est bien meilleur. Autrement dit, si je cherche à approximer un filtre passe-bande pas trop serré autour de n'importe quelle valeur propre, l'approximation est bonne (et ne dépend plus de la densité spectrale). Aussi, le paragraphe "notations" est à adapter en fonction des notations dont on a eu effectivement besoin pour l'article.}

\section{Acknowledgments}
This  work  has  been  partially  supported  by  the  LabEx PERSYVAL-Lab  (ANR-11-LABX-0025-01)  funded  by  the French  program  Investissement  d'avenir.

\bibliographystyle{IEEEtran}
\bibliography{IEEEabrv,biblio}

\begin{thebibliography}{10}
\providecommand{\url}[1]{#1}
\csname url@samestyle\endcsname
\providecommand{\newblock}{\relax}
\providecommand{\bibinfo}[2]{#2}
\providecommand{\BIBentrySTDinterwordspacing}{\spaceskip=0pt\relax}
\providecommand{\BIBentryALTinterwordstretchfactor}{4}
\providecommand{\BIBentryALTinterwordspacing}{\spaceskip=\fontdimen2\font plus
\BIBentryALTinterwordstretchfactor\fontdimen3\font minus
  \fontdimen4\font\relax}
\providecommand{\BIBforeignlanguage}[2]{{%
\expandafter\ifx\csname l@#1\endcsname\relax
\typeout{** WARNING: IEEEtran.bst: No hyphenation pattern has been}%
\typeout{** loaded for the language `#1'. Using the pattern for}%
\typeout{** the default language instead.}%
\else
\language=\csname l@#1\endcsname
\fi
#2}}
\providecommand{\BIBdecl}{\relax}
\BIBdecl

\bibitem{Shuman2013}
D.~I. Shuman, S.~K. Narang, P.~Frossard, A.~Ortega, and P.~Vandergheynst, ``The
  emerging field of signal processing on graphs: Extending high-dimensional
  data analysis to networks and other irregular domains,'' \emph{Signal
  Processing Magazine, IEEE}, vol.~30, no.~3, 2013.

\bibitem{Sandryhaila2013}
A.~Sandryhaila and J.~Moura, ``Discrete signal processing on graphs,''
  \emph{Signal Processing, IEEE Transactions on}, vol.~61, no.~7, 2013.

\bibitem{tremblay_gsp_chapter}
\BIBentryALTinterwordspacing
N.~Tremblay, P.~Goncalves, and P.~Borgnat, ``Design of graph filters and
  filterbanks,'' oct 2017. [Online]. Available:
  \url{https://hal.archives-ouvertes.fr/hal-01243889}
\BIBentrySTDinterwordspacing

\bibitem{CooleyTukey1965}
J.~Cooley and J.~Tukey, ``An algorithm for the machine calculation of complex
  {F}ourier series,'' \emph{Mathematics of Computation}, vol.~19, no.~90, pp.
  297--301, 1965.

\bibitem{Morgenstern1975}
J.~Morgenstern, ``The linear complexity of computation,'' \emph{J. ACM},
  vol.~22, no.~2, pp. 184--194, Apr. 1975.

\bibitem{Lemagoarou2016a}
L.~Le~Magoarou and R.~Gribonval, ``Flexible multilayer sparse approximations of
  matrices and applications,'' \emph{IEEE Journal of Selected Topics in Signal
  Processing}, vol.~10, no.~4, pp. 688--700, 2016.

\bibitem{magoarou_approximate_2017}
L.~L. Magoarou, R.~Gribonval, and N.~Tremblay, ``Approximate fast graph
  {Fourier} transforms via multi-layer sparse approximations,'' \emph{IEEE
  Transactions on Signal and Information Processing over Networks}, vol.~PP,
  no.~99, pp. 1--1, 2017.

\bibitem{Givens1958}
W.~Givens, ``\BIBforeignlanguage{English}{Computation of plane unitary
  rotations transforming a general matrix to triangular form},''
  \emph{\BIBforeignlanguage{English}{Journal of the Society for Industrial and
  Applied Mathematics}}, vol.~6, no.~1, pp. pp. 26--50, 1958.

\bibitem{Jacobi1846}
C.~G.~J. Jacobi, ``\"{U}ber ein leichtes verfahren, die in der theorie der
  s\"{a}kularst\"{o}rungen vorkommenden gleichungen numerisch aufzul\"{o}sen,''
  \emph{J. reine angew. Math.}, vol.~30, pp. 51--94, 1846.

\bibitem{Golub2000}
G.~H. Golub and H.~A. Van~der Vorst, ``Eigenvalue computation in the 20th
  century,'' \emph{Journal of Computational and Applied Mathematics}, vol. 123,
  no.~1, pp. 35--65, 2000.

\bibitem{Golub2012}
G.~H. Golub and C.~F. Van~Loan, \emph{Matrix computations}.\hskip 1em plus
  0.5em minus 0.4em\relax JHU Press, 2012, vol.~3.

\bibitem{decelle_asymptotic_2011}
A.~Decelle, F.~Krzakala, C.~Moore, and L.~Zdeborová, ``Asymptotic analysis of
  the stochastic block model for modular networks and its algorithmic
  applications,'' \emph{Phys. Rev. E}, vol.~84, no.~6, p. 066106, 2011.

\end{thebibliography}
\end{document}